\documentclass{article}
\usepackage{amsmath, amsthm}
\usepackage{amssymb}
\newtheoremstyle{theorem}
 {10pt}          
 {10pt}  
 {\sl}  
 {\parindent}     
 {\bf}  
 {. }    
 { }    
 {}     
\theoremstyle{theorem}
\newtheorem{theorem}{Theorem}

\newtheorem{lemma}[theorem]{Lemma}
\newtheorem{remark}[theorem]{Remark}
\newtheorem{proposition}[theorem]{Proposition}
\newtheorem{example}[theorem]{Example}
\newtheoremstyle{defi}
 {10pt}          
 {10pt}  
 {\rm}  
 {\parindent}     
 {\bf}  
 {. }    
 { }    
 {}     
\theoremstyle{defi}
\newtheorem{definition}[theorem]{Definition}

\begin{document}

\title{A Quantum  Generalized  Mittag-Leffler Function Via Caputo
q-Fractional  Equations}

\author{Thabet Abdeljawad and Bet\"{u}l Benli\\
 Department of Mathematics and
Computer Science\\
 \c{C}ankaya University, 06530 Ankara, Turkey}
\maketitle
\begin{abstract}
Some Caputo q-fractional difference equations  are solved. The solutions
are expressed  by means of  a new introduced generalized type of
q-Mittag-Leffler  functions. The method of successive approximation is
used to obtain the solutions. The obtained q-version of Mittag-Leffler
function is thought as the q-analogue of the one introduced previously by
Kilbas and Saigo.



\end{abstract}

\section{Introduction and Preliminaries} \label{s:1}
The concept of fractional calculus is not new.  However, it has gained its
popularity and importance during the last three decades or so. This is due
to its distinguished applications in numerous diverse fields of science
and engineering (\cite{Samko}, \cite{Podlubny}, \cite{Kilbas}). The
q-calculus is also not of recent appearance. It was initiated in twenties
of the last century.For the basic concepts in q- calculus we refer the
reader to \cite{history}. Starting from the q-analogue of Cauchy formula
\cite{Al-salam2}, Al-Salam started the fitting of the concept of
q-fractional calculus. After that he (\cite{Alsalaml}, \cite{Alsalam}) and
Agarwal R. \cite{Agarwal} continued on by studying certain q-fractional
integrals and derivatives, where they proved the semigroup properties for
left and right (Riemann)type fractional integrals but without variable
lower limit and variable upper limit, respectively. Recently, the authors
in \cite{Pred} generalized the notion of the (left)fractional q-integral
and q-derivative by introducing variable lower limit and proved the
semigroup properties.

Very recently and after the appearance of time scale calculus (see for
example \cite{Boh}), some authors started to pay attention and apply the
techniques of time scale to discrete fractional calculus
(\cite{Ferd},\cite{Feri},\cite{Ferq}, \cite{Th}) benefitting from the
results announced before in \cite{Miller}. All of these results are mainly
about fractional calculus on the time scales $T_q=\{q^n:n \in
\mathbb{Z}\}\cup \{0\}$ and $h\mathbb{Z}$ \cite{Bastos}. As a contribution
in  this direction and being motivated by all above, in this article we
introduce the q-analogue of a generalized type Mittag-Leffler function
used before by Kilbas and Saigo in \cite{Saigo}. Such  functions are
obtained by solving  linear q-Caputo initial value problems. The results
obtained in this article generalize also the results of \cite{TD2011}.

For the theory of q-calculus we refer the reader to the survey
\cite{history} and for the basic definitions and results for the
q-fractional calculus we refer to \cite{Ferq}. Here we shall
summarize some of those basics.

For $0<q<1$, let $T_q$ be the time scale
$$T_q=\{q^n:n \in \mathbb{Z}\}\cup \{0\}.$$ where $Z$ is the set of
integers. More generally, if $\alpha$ is a nonnegative real number
then we define the time scale
$$T_q^\alpha=\{q^{n+\alpha}:n \in Z\}\cup \{0\},$$ we write $T_q^0=T_q.$

 For a function $f:T_q\rightarrow \mathbb{R}$, the nabla
q-derivative of $f$ is given by
\begin{equation} \label{qd}
\nabla_q f(t)=\frac{f(t)-f(qt)}{(1-q)t},~~t \in T_q-\{0\}
\end{equation}
The nabla q-integral of $f$ is given by
\begin{equation} \label{qi}
\int_0^t f(s)\nabla_q s=(1-q)t\sum_{i=0}^\infty q^if(tq^i)
\end{equation}
and for $0\leq a \in T_q$

$$\int_a^t f(s)\nabla_q s=\int_0^t f(s)\nabla_q s - \int_0^a f(s)\nabla_q s$$
On the other hand
\begin{equation} \label{r1}
\int_t^\infty f(s)\nabla_q s=(1-q)t\sum_{i=1}^\infty q^{-i}
f(tq^{-i})
\end{equation}
 and for $0<b<\infty$ in $T_q$
\begin{equation} \label{r2}
\int_t^b f(s)\nabla_q s=\int_t^\infty f(s)\nabla_q s - \int_b^\infty
f(s)\nabla_q s
\end{equation}
 By the fundamental theorem in q-calculus we have
\begin{equation} \label{fq}
\nabla_q \int_0^t f(s)\nabla_q s=f(t)
\end{equation}
and if $f$ is continuous at $0$, then
\begin{equation} \label{cfq}
\int_0^t \nabla_qf(s)\nabla_q s=f(t)-f(0)
\end{equation}
Also the following identity will be helpful
\begin{equation} \label{help}
\nabla_q\int_a^t f(t,s)\nabla_q s=\int_a^t \nabla_qf(t,s)\nabla_q
s+f(qt,t)
\end{equation}
Similarly the following identity will be useful as well
\begin{equation} \label{help1}
\nabla_q\int_t^b f(t,s)\nabla_q s=\int_{qt}^b \nabla_qf(t,s)\nabla_q
s-f(t,t)
\end{equation}
The q-derivative in (\ref{help}) and (\ref{help1}) is applied with
respect to t.

From the theory of q-calculus and the theory of time scale more
generally, the following product rule is valid

\begin{equation} \label{qproduct}
\nabla_q (f(t)g(t))=f(qt)\nabla_q g(t)+\nabla_q f(t)g(t)
\end{equation}
The q-factorial function for $n\in \mathbb{N}$ is defined by
\begin{equation} \label{qfact}
(t-s)_q^n=\prod_{i=0}^{n-1}(t-q^is)
\end{equation}
When $\alpha$ is a non positive integer, the  q-factorial function
is defined by
\begin{equation} \label{qfactg}
(t-s)_q^\alpha=t^\alpha\prod_{i=0}^\infty \frac{1- \frac{s}{t} q^i}
{1- \frac{s}{t} q^{i+\alpha}}
\end{equation}
We summarize some of the properties of  q-factorial functions, which
can be found mainly in \cite{Ferq}, in the following lemma

\begin{lemma} \label{qproperties}
(i)$ (t-s)_q^{\beta+\gamma}=(t-s)_q^\beta (t-q^\beta s)_q^\gamma$

(ii)$(at-as)_q^\beta=a^\beta (t-s)_q^\beta$

(iii) The nabla q-derivative of the q-factorial function with
respect to $t$ is

$$\nabla_q (t-s)_q^\alpha
=\frac{1-q^\alpha}{1-q}(t-s)_q^{\alpha-1}$$

(iv)The nabla q-derivative of the q-factorial function with respect
to $s$ is
$$\nabla_q (t-s)_q^\alpha
=-\frac{1-q^\alpha}{1-q}(t-qs)_q^{\alpha-1}$$ where
$\alpha,\gamma,\beta \in \mathbb{R}.$
\end{lemma}

\begin{definition}\label{qcd}\cite{TD2011}
Let  $\alpha>0$. If $\alpha \notin \mathbb{N}$ , then the
$\alpha-$order Caputo (left) q-fractional
derivative of a function $f$ is defined by

\begin{equation} \label{qrd}
_{q}C_a^\alpha f(t)\triangleq~ _{q}I_a ^{(n-\alpha)}\nabla_q
^nf(t)=\frac{1}{\Gamma(n-\alpha)} \int_a^t(t-qs)_q^{n-\alpha-1}
\nabla_q^nf(s)\nabla_q s
\end{equation}

where $n=[\alpha]+1$.

If $\alpha \in \mathbb{N}$, then $_{q}C_a^\alpha f(t)\triangleq
\nabla_q^n f(t)$
\end{definition}
It is clear that $_{q}C_a^\alpha $ maps functions defined on
$T_q$ to functions defined on $T_q$, and that $_{b}C_q^\alpha $ maps
functions defined on $T_q^{1-\alpha}$ to functions defined on $T_q$

The following identity which  is useful to transform Caputo q-fractional
difference equations  into q-fractional integrals, will be  our key in
solving the q-fractional linear type equation by using successive
approximation.

\begin{proposition} \label{qtranss}\cite{TD2011}
Assume $\alpha>0$ and $f$ is defined in suitable domains. Then
\begin{equation}\label{qtrans1}
_{q}I_a^\alpha~ _{q}C_a^\alpha
f(t)=f(t)-\sum_{k=0}^{n-1}\frac{(t-a)_q^{k}}{\Gamma_q(k+1)}\nabla_q^kf(a)
\end{equation}
and if $0<\alpha\leq1$ then
\begin{equation}\label{qtrans3}
_{q}I_a^\alpha~ _{q}C_a^\alpha f(t)= f(t)-f(a)
 \end{equation}

\end{proposition}

The following identity \cite{Pred} is essential to solve linear
q-fractional equations
\begin{equation}\label{qpower}
_{q}I_a^\alpha (x-a)_q^\mu
=\frac{\Gamma_q(\mu+1)}{\Gamma_q(\alpha+\mu+1)}(x-a)_q^{\mu+\alpha}~~(0<a<x<b)
\end{equation}
where $\alpha \in \mathbb{R}^+$ and $\mu \in (-1,\infty)$.
The q-analogue of Mittag-Leffler function with double index
$(\alpha,\beta)$ is introduced in \cite{TD2011}. It was defined as
follows:

\begin{definition}\label{Mitt}\cite{TD2011}
For $z ,z_0\in \mathbf{C}$ and $\mathfrak{R}(\alpha)> 0$, the
q-Mittag-Leffler function is defined by
\begin{equation} \label{Mit}
_{q}E_{\alpha,\beta}(\lambda,z-z_0)=\sum_{k=0}^\infty \lambda^k
\frac{(z-z_0)_q^{\alpha k}}{\Gamma_q(\alpha k+\beta)}.
\end{equation}
When $\beta=1$ we simply use
$~_{q}E_{\alpha}(\lambda,z-z_0):=~_{q}E_{\alpha,1}(\lambda,z-z_0)$.
\end{definition}

\section{Main Results } \label{s:1}
The following is to be the q-analogue of the generalized Mittag-Leffler
function introduced by Kilbas and Saigo \cite{Saigo} (see also
\cite{Kilbas} page 48).
\begin{definition}
For $\alpha,l,\lambda \in \mathbb{C}$ are complex numbers and $m \in
\mathbb{R}$ such that $\mathfrak{R}(\alpha)>0, ~m >0,a\geq 0 $ and
$\alpha(jm+l)\neq -1,-2,-3,...,~$, the generalized q-Mittag-Leffler
function (of order 0) is defined by
$$_{q}E_{\alpha,m,l}(\lambda,x-a)=1+ \sum_{k=1}^\infty \lambda^k
q^{-\frac{k(k-1)}{2}  \alpha (m-1)(\alpha l+\alpha)}c_k (x-a)_q^{\alpha k
m}$$ where
 $$c_k=\prod_{j=0}^{k-1}\frac{\Gamma_q[\alpha(jm+l)+1]}{\Gamma_q[\alpha(jm+l+1)+1]},~k=1,2,3,...$$
 While the the generalized q-Mittag-Leffler function (of order r),
$r=0,1,2,3,...$, is defined by
 $$_{q}E^r_{\alpha,m,l}(\lambda,x-a)=1+ \sum_{k=1}^\infty \lambda^k q^{-k
\alpha (m-1)r} q^{-\frac{k(k-1)}{2}  \alpha (m-1)(\alpha l+\alpha)}c_k
(x-q^r a)_q^{\alpha k m}.$$
 Note that
$_{q}E^0_{\alpha,m,l}(\lambda,x-a)=~_{q}E_{\alpha,m,l}(\lambda,x-a).$
\end{definition}
\begin{remark}\label{r}
In particular, if $m=1$, then the generalized q-Mittag-Leffler function is
reduced to the q-Mittag-Leffler function, apart from a constant factor
$\Gamma_q(\alpha l+1)$. Namely,
\begin{equation} \label{sp}
_{q}E_{\alpha,1,l}(\lambda,x-a)=\Gamma_q(\alpha l+1)~_{q}E_{\alpha,\alpha
l+1}(\lambda,x-a)
 \end{equation}

 This turns to be the q-analogue of the identity
$E_{\alpha,1,l}(z)=\Gamma(\alpha l+1)E_{\alpha,\alpha l+1}(z)$ (see
\cite{Kilbas}) page 48).
\end{remark}

\begin{example} \label{q}
Consider the q-Caputo difference equation
\begin{equation} \label{mainq}
(_{q}C_a^\alpha y)(x)= \lambda (x-a)_q^\beta y(q^{-\beta}x),~y(a)=b
\end{equation}
where
$$~~~0< \alpha <1,~\beta>-\alpha,~\lambda \in \mathbb{R},~b \in \mathbb{R}.$$

Applying Proposition \ref{qtranss} we have
$$y(x)=y(a)+ \lambda~ _{q}I_a^\alpha [(x-a)_q^\beta y(q^{-\beta}x)].$$

The method of successive applications implies that
$$y_m (x)= y(a)+ \lambda _{q}I_a^\alpha [(x-a)_q^\beta
y_{m-1}(q^{-\beta}x)],~~m=1,2,3,...,$$
where $y_0(x)=b$. Then by the help of (\ref{qpower}) we have
$$y_1(x)=b+b \lambda
\frac{\Gamma_q(\beta+1)}{\Gamma_q(\beta+\alpha+1)}(x-a)_q^{\beta+\alpha},$$
and
$$y_2(x)=b+b \lambda~ _{q}I_a^\alpha [(x-a)_q^\beta\{1+ \lambda
\frac{\Gamma_q(\beta+1)}{\Gamma_q(\beta+\alpha+1)}(q^{-\beta}x-a)_q^{\beta+\alpha}
\}] $$
Then by (i) and (ii)  of Lemma \ref{qproperties}
$$y_2(x)=b+b \lambda~ _{q}I_a^\alpha [(x-a)_q^\beta+ \lambda
\frac{\Gamma_q(\beta+1)}{\Gamma_q(\beta+\alpha+1)}
q^{-\beta(\alpha+\beta)}(x-a)_q^{2\beta+\alpha}]$$
Again by (\ref{qpower}) we conclude
$$y_2(x)=b+b \lambda~ _{q}I_a^\alpha [(x-a)_q^\beta+ \lambda
\frac{\Gamma_q(\beta+1)}{\Gamma_q(\beta+\alpha+1)}
q^{-\beta(\alpha+\beta)}(x-a)_q^{2\beta+\alpha}]$$
Then (\ref{qpower}) leads to
$$y_2(x)=$$
\begin{equation}
b[1+\lambda~ \frac{\Gamma_q(\beta+1)}{\Gamma_q^(\beta+\alpha+1)}
(x-a)_q^{\beta+\alpha}+ \lambda^2
\frac{\Gamma_q(2\beta+\alpha+1)}{\Gamma_q(2\beta+2\alpha+1)}
q^{-\beta(\alpha+\beta)}(x-a)_q^{2\beta+2\alpha}].
\end{equation}
Proceeding inductively, for each $m=1,2,..$ we obtain

\begin{equation}
y_m(x)=b \large[1+ \sum_{k=1}^m \lambda^k
q^{-\beta\frac{k(k-1)}{2}(\alpha+\beta)}
c_k(x-a)_q^{k(\alpha+\beta)}\large]
\end{equation}
where
$$c_k=\prod_{j=0}^{k-1}\frac{\Gamma_q[\alpha(jm+l)+1]}{\Gamma_q[\alpha(jm+l+1)+1]},~~m=1+\frac{\beta}{\alpha},~l=\frac{\beta}{\alpha},~k=1,2,3,...$$
\end{example}
If we let $m\rightarrow \infty$, then we obtain the solution
$$y(x)=b ~\large[1+ \sum_{k=1}^\infty \lambda^k
q^{-\beta\frac{k(k-1)}{2}(\alpha+\beta)}
c_k(x-a)_q^{k(\alpha+\beta)}\large]$$
which is exactly
$$y(x)=b~
_{q}E_{\alpha,1+\frac{\beta}{\alpha},\frac{\beta}{\alpha}}(\lambda,x-a).$$

\begin{remark} \label{rr}
1)If in (\ref{mainq}) $\beta=0$, then in accordance with (\ref{sp}) and
Example 10 in \cite{TD2011} we have
$$_{q}E_{\alpha,1,0}(\lambda,x-a)=~_{q}E_{\alpha,1}(\lambda,x-a)=~_{q}E_{\alpha}(\lambda,x-a)$$

2) The solution of the q-Cauchy problem
\begin{equation} \label{mainqhalf}
(_{q}C_a^{\frac{1}{2}} y)(x)= \lambda (x-a)_q^\beta y(q^{-\beta}x),~y(a)=b
\end{equation}
where
$$~~~0< \alpha <1,~\beta>-\frac{1}{2},~\lambda \in \mathbb{R},~b \in
\mathbb{R}$$ is given by
$$y(x)=b~ _{q}E_{\frac{1}{2},1+2\beta,2\beta}(\lambda,x-a).$$

3) By the help of (\ref{qtrans1}) and Lemma \ref{qproperties} and by
applying the successive approximation with
$y_0(x)=\sum_{k=0}^{n-1}\frac{(t-a)_q^{k}}{\Gamma_q(k+1)}\nabla_q^kf(a)$,
Example \ref{q} can be generalized for arbitrary $\alpha >0$. Namely, the
solution of the q-initial value problem

\begin{equation} \label{gmainq}
(_{q}C_a^\alpha y)(x)= \lambda (x-a)_q^\beta
y(q^{-\beta}x),~y^{(k)}(a)=b_k ~~(b_k \in \mathbb{R},~k=0,1,...,n-1)
\end{equation}
where
$$~~~n-1< \alpha < n ,~\beta>-\alpha,~\lambda \in \mathbb{R},~b \in
\mathbb{R}$$ is given by

$$y(x)= \sum_{r=0}^{n-1}\frac{b_r}{\Gamma_q(r+1)}(x-a)_q^r~
_{q}E^r_{\alpha,1+\frac{\beta}{\alpha},\frac{\beta
+r}{\alpha}}(\lambda,x-a).$$ Note that when $0< \alpha < 1$, i.e, $n=1$,
the solution of Example \ref{q} is recovered.
\end{remark}

\end{document}